\documentclass[preprint,10pt]{elsarticle}
\usepackage{amssymb}
\usepackage{amsthm,amsmath}
\usepackage{filecontents, color, url}
\usepackage{natbib}

\newtheorem{thm}{Theorem}[section]
\newtheorem{prop}[thm]{Proposition}
\newtheorem{lem}[thm]{Lemma}

\newtheorem{defn}[thm]{Definition}

\newcommand{\R}{\mathbb{R}}  
\newcommand{\N}{\mathbb{N}}  
\newcommand{\sw}{\mathcal{S}} 
\newcommand{\domain}{[0,2\pi]\times \R} 

\journal{}

\begin{document}

\begin{frontmatter}


\title{\large{Modified Radon transform inversion using moments in $\R^2$}}

\author[label1]{Hayoung Choi}
\ead{hayoung79choi@gmail.com}
\author[label2]{Farhad Jafari\corref{cor1}}
\ead{fJafari@uwyo.edu}
\author[label3]{Robert Mnatsakanov}
\ead{rmnatsak@stat.wvu.edu\fnref{label3}}

\address[label1]{School of Information Science and Technology, ShanghaiTech University, Shanghai 201210, China}
\address[label2]{Department of Mathematics, University of Wyoming, Laramie, WY 82070, USA}
\address[label3]{Department of Statistics, West Virginia University
Morgantown, WV 26506, USA}


\begin{abstract}
A modified Radon transform for noisy data is introduced and its inversion formula is established.
The problem of recovering the multivariate probability density function $f$ from the moments of its modified Radon transform $\widehat{R}f$ is considered. 
\end{abstract}


\end{frontmatter}



\section{Introducton}

Radon transform of an integrable function is the integral of that function over lines. A key application of Radon transform is tomography where the interior density of a 2-D object (e.g. slices of a 3-D object) is reconstructed from the Radon transform data. In this paper, we consider the problem of recovering a non-negative desnity function from noisy Radon transform data using moments of this data.

In general, we will follow the notation in \cite{book:Helgason} and \cite{Quinto06}. The Radon transform $R$ in $\R^2$ is defined by
\begin{equation}
Rf(w,p)=\int_{\langle x,\omega\rangle=p} f(x)dm(x),\label{defRadon1}
\end{equation}
where $\omega=(\omega_1,\omega_2)$ is a unit vector, $p\in \R$, and $dm$ is the arc length measure on the line $\langle x,\omega\rangle =p$ with the usual inner product $\langle~,~\rangle$.
Clearly, the Radon transform can be represented as an integral transform with respect to a measure $ \mu $, which is singular with respect of the Lebesgue measure in $ \R^2 $.
\begin{equation}
Rf(w,p)=\int_{\R^n}  f(x) d \mu = \int_{\R^n} f(x) \chi_{\{\langle x,\omega\rangle =p\}} dx,\label{defRadon2}
\end{equation}
The measure $ \mu $ restricts the Lebesgue measure to lines $ E $, parameterized by $ \omega $ and $ p $, in $ \R^2 $ and $\chi_E$ denotes the indicator function of the set $ E $.
A function $f$ is said to be in the \emph{Schwartz space} $\mathcal{S}(\R^2)$ if $f\in C^\infty(\R^2)$ and for each integer $ m \geq 0 $ and each polynomial $P$ of homogeneous degree $ m $
\begin{equation*}
\underset{x}{\textup{sup}} \big||x|^\alpha P(\partial_{x_1}, \partial_{x_2})f(x)\big|<\infty,
\end{equation*}
where $|x|$ is the Euclidean norm of $x$, and $ \alpha $ is a multi-index.
 A function $g(\theta,p)$ is said to be in the Schwartz space $\mathcal{S}([0, 2\pi] \times \R)$ if $g(\theta,p)$ can be extended to a smooth and $2\pi$--periodic function in $\theta$, and $g(\cdot,p)\in \mathcal{S}(\R)$ uniformly in $\theta$.  
 
In the remainder of this section, we will collect a few well known results about the Radon transform and the Hamburger moment problem that will be needed in the remainder of this paper.  For the sake of completeness, these results are stated. Readers are referred to references for deeper treatment of these results.

\begin{lem} \textup{(see \cite{book:Helgason}, Lemma 2.3)}. \label{range} 
For each $f\in \mathcal{S}(\R^2)$ the Radon transform $Rf$ satisfies the following condition: For $k\in \N_0$ the integral
\begin{equation*}
\int_{\R} Rf(\omega,p)p^k~dp
\end{equation*}
is a $k^{th}$ degree homogeneous polynomial in $\omega_1,\omega_2$.
\end{lem}

We denote the unit vector in direction $\theta$ as $\omega=\omega(\theta):=(\omega_1,\omega_2)$ with $\omega_1=\cos{\theta}$ and  $\omega_2=\sin{\theta}$. Thus, the Radon transform of $f\in L^1(\R^2)$ can be expressed as a function of $(\theta,p)$:
\begin{equation}
Rf(\theta,p)=\int_{\langle x,\omega(\theta)\rangle=p} f(x)dx.\label{defRadon3}
\end{equation}

Note that since the pairs $(\omega,p)$ and $(-\omega,-p)$ give the same line, $R$ satisfies the \emph{evenness} condition: $Rf(\theta,p)=Rf(\theta+\pi,-p)$.

\begin{thm}\label{radon:cont}
The Radon transform $R$ is a bounded linear operator from $L^1(\R^2)$ to $L^1([0,2\pi]\times \R)$
with norm $\|R\|\leq 2\pi$,
i.e., $\|R f\|_{L^1([0,2\pi]\times \R)} \leq 2\pi \|f\|_{L^1(\R^2)}$.
\begin{proof}
See \cite{Quinto06}, for example.
\end{proof}
\end{thm}

Along with the transform $Rf$ we define the dual Radon transform of $g\in L^1([0,2\pi] \times \R)$ as
\begin{equation}
R^*g(x)=\int_{0}^{2\pi} g(\theta, \langle x, \omega \rangle) d\theta,
\end{equation}
which is the integral of $g$ over all lines through $x$.

Using $F_1$ and $F_2$ for the 1-D and 2-D Fourier transforms, recall that
\begin{align*}
F_1f(s) &= \frac{1}{\sqrt{2\pi}} \int_{-\infty}^{\infty} f(t) e^{-ist} dt,\\
F_1^{-1} f(t) &= \frac{1}{\sqrt{2\pi}} \int_{-\infty}^{\infty} f(s) e^{ist} ds,\\
F_2f(\xi) &= \frac{1}{2\pi} \int_{\R^2} f(x) e^{-i\langle x,\xi\rangle } dx,\\
F_2^{-1}f(x) &= \frac{1}{2\pi} \int_{\R^2} f(\xi) e^{i\langle \xi,x\rangle } d\xi.
\end{align*}

\begin{thm}[Projection-Slice Theorem]\label{thm:pojection-slice}
Let $f\in L^1(\R^2)$.
Then,
\begin{equation*}
F_2 f(s\omega)=\frac{1}{\sqrt{2\pi}} F_1 (Rf(\theta,\cdot))(s).
\end{equation*}
\end{thm}
This theorem shows that $R$ is injective on $L^1(\R^2)$.
The Fourier inversion formula combined with the Projection-Slice Theorem provides an inversion formula for the Radon transform in $ \R^2 $.

Denote the Riesz potential $I^{-1}$, for $g\in L^1([0,2\pi]\times \R)$, as the operator with Fourier multiplier
$|s|$ (see \cite{book:Stein}):
\begin{equation}
I^{-1}g=F_1^{-1}(|s|(F_1g(\theta,\cdot)(s))).
\end{equation}

\begin{thm}[Inversion formula for $R f$]\label{thm:inversionf}
Let $f\in C_c^\infty(\R^2)$. Then
\begin{equation*}
f(x)=\frac{1}{4\pi}R^* (I^{-1} Rf)(x).
\end{equation*}
\end{thm}
Note that this theorem is true on a larger domain than $C_c^\infty(\R^2)$. However, $I^{-1}Rf$ may be a distribution rather than a function.
%

One may ask, when is a given function $g$ the Radon transform of a function $f$? In other words, for a given function $g$, does there exist $f$ such that $g=Rf$? The following theorem is the fundamental result on this question, which is called the \emph{Schwartz theorem}
or \emph{Range Theorem} for the Radon transform.

\begin{thm}
Let $g\in \sw(\domain)$ be even. Then, there exists $f\in \sw(\R^2)$ such that $g=Rf$
if and only if
\begin{align}
& \textup{for each } k\in \N_0,  \textup{ the $k^{th}$ moment }
\int_{-\infty}^{\infty} g(\theta,p)p^kdp \label{homo}\\
& \textup{ is a homogeneous polynomial of degree $k$ in } \omega_1 \textup{ and } \omega_2.\nonumber
\end{align}
\begin{proof}
See \cite{book:Helgason}, Theorem 2.4.
\end{proof}
\end{thm}

In many applications, the collected dat is not the original density function $ f$ but $ Rf $ or, in fact $ Rf + \eta $, where $ \eta $ is spurious noise attributable to non-zero width of the pencil beam, motion of the object being interrogated and detector noise. We propose a method to approximation the density function $ f $ from its noisy Radon transform using moments as follows:\\
1. Collect $Rf$ including noisy data;\\
2. Find the modified Radon transform $\widehat{R}f$;\\
3. Recover moments from the modified Radon transform;\\
4. Recover $f$ via its moments. \\
This paper will establish the details of this reconstruction.

%
%
\section{Modified Radon Transform}

\begin{defn}\label{defn:molifier}
If $\varphi$ is a smooth function on $\R$, satisfying the following requirements:\\
(1) it is compactly supported,\\
(2) ${\displaystyle \int_{-\infty}^{\infty} \varphi(t)dt=1}$,\\
(3) ${\displaystyle  \lim_{\varepsilon\rightarrow 0} \varphi _{\varepsilon}(t) :=
 \lim_{\varepsilon\rightarrow 0} \frac{1}{\varepsilon} \varphi (\frac{t}{\varepsilon})=\delta(t) }$,\\
(4) ${\displaystyle \varphi(t)\geq 0 }$ for all $t\in\R $,\\
then $ \varphi $ is called a \textit{positive mollifier}.\\
Furthermore, if \\
(5) $\varphi(t) = h (|t|)$ for some infinitely differentiable function
$h :\R^n\rightarrow \R$, then it is called a \textit{symmetric mollifier}.
\end{defn}

For example, if $\varphi:\R\rightarrow \R$ is the Gaussian function
\begin{equation*}
\varphi(x) = \left\{
\begin{array}{rl}
e^{-\frac{1}{1-|x|^2}} & \text{if } |x| < 1,\\
0 & \text{if } |x| \geq 1.
\end{array} \right.
\end{equation*}
then $\varphi$ is a positive symmetric mollifier.

Let $\Omega=\{\varphi : \varphi$ is a positive symmetric mollifier such that
$F_1(\varphi)(s)>0$ for all $s\in\R \}$. Note that $F_1$ is the Fourier transform.
Clearly, the Gaussian function is in $\Omega$.

\begin{defn}
Let $\varphi\in \Omega$ and $f\in L^1(\R^2)$.
The \emph{modified Radon transform} in $\R^2$ is defined by
\begin{equation}
\widehat{R}_\varphi f(\theta,p) = \int_{\R^2} \chi(x;\theta,p) f(x)dx,
\end{equation}
where
\begin{equation*}
\chi(x;\theta,p) = (\delta \ast \varphi) (\langle x,\omega \rangle-p),~ (\delta \textup{ is a delta function }).
\end{equation*}
\end{defn}

It is easy to show that modified Radon transform, which is defined as the Radon transform of the smoothed density, is identical to the smoothed Radon transform of the original density.  

\begin{lem}\label{modified:form1}
\begin{equation}
\widehat{R}_\varphi f(\theta,p) = \int_{\R} Rf(\theta,p+\tau) \varphi (\tau) d\tau
\end{equation}
\begin{proof}
\begin{align*}
\widehat{R}_\varphi f(\theta,p)
& = \int_{\R^2} (\delta \ast \varphi)(\langle x,\omega\rangle-p) f(x)dx\\
& = \int_{\R} \Bigg( \int_{\R^2} \delta(\langle x,\omega \rangle-p-\tau) f(x)dx \Bigg) \varphi (\tau) d\tau \\
& = \int_{\R} \Bigg( \int_{\langle x,\omega \rangle=p+\tau} f(x)dx \Bigg) \varphi (\tau) d\tau.
\end{align*}
\end{proof}
\end{lem}

\begin{prop}\label{modified:form2}
\begin{equation*}
\widehat{R}_\varphi f(\theta,p) = \Big(Rf(\theta,\cdot) \ast \varphi \Big)(p).
\end{equation*}
\begin{proof}
By Lemma \ref{modified:form1}
\begin{align*}
\widehat{R}_\varphi f(\theta,p)
& = \int_{\R} Rf(\theta,p+\tau) \varphi (\tau) d\tau \\
& = \int_{\R} Rf(\theta,p-\tau) \varphi (-\tau) d\tau.
\end{align*}
Since $\varphi$ is a symmetric mollifier,
\begin{equation*}
\widehat{R}_\varphi f(\theta,p)
= \int_{\R} Rf(\theta,p-\tau) \varphi (\tau) d\tau.
\end{equation*}
\end{proof}
\end{prop}

\begin{thm}\label{boundedL1}
The modified Radon transform $\widehat{R}_\varphi$ is a bounded linear operator from $L^1(\R^2)$ to $L^1([0,2\pi]\times \R)$ with norm $\|\widehat{R}_\varphi\|\leq 2\pi$,i.e.,
$\|\widehat{R}_\varphi f\|_{L^1([0,2\pi]\times \R)} \leq 2\pi \|f\|_{L^1(\R^2)}$.
\begin{proof}
Using Proposition \ref{modified:form2}, one finds that
\begin{align*}
\|\widehat{R}_\varphi f\|_{L^1([0,2\pi] \times \R)}
&= \int_{\theta=0}^{2\pi}  \int_{p=-\infty}^{\infty}  \Big|\widehat{R}_\varphi f(\theta,p)\Big|dp d\theta\\
&= \int_{\theta=0}^{2\pi}  \int_{p=-\infty}^{\infty}
\Bigg|\int_{\tau=-\infty}^{\infty} Rf(\theta, p+\tau) \varphi (\tau) d\tau \Bigg| dp d\theta\\
&\leq \int_{\theta=0}^{2\pi}  \int_{p=-\infty}^{\infty}
\int_{\tau=-\infty}^{\infty} \Big|Rf(\theta, p+\tau)\Big| \varphi (\tau) d\tau  dp d\theta\\
&= \int_{\tau=-\infty}^{\infty} \varphi(\tau) \Bigg( \int_{\theta=0}^{2\pi}  \int_{p=-\infty}^{\infty}
 \Big|Rf(\theta, p+\tau)\Big| dp d\theta \Bigg)  d\tau.
\end{align*}
By Theorem \ref{radon:cont} and Definition \ref{defn:molifier}, it follows that
\begin{align*}
\|\widehat{R}_\varphi f\|_{L^1([0,2\pi] \times \R)}
&\leq \int_{\tau=-\infty}^{\infty} \varphi(\tau) \Big( 2\pi \|f\|_{L^1(\R^2)} \Big)  d\tau\\
&= 2\pi \|f\|_{L^1(\R^2)}.
\end{align*}
\end{proof}
\end{thm}

Denote the modified Riesz potential $\widehat{I}^{-1}$, for $g\in L^1([0,2\pi]\times \R)$, as the operator with Fourier multiplier $|s|$ and symmetric mollifier $\varphi$:
\begin{equation*}
\widehat{I}^{-1}g=F_1^{-1}\bigg(|s|\Big(\frac{F_1(g(\theta,\cdot))(s)}{F_1(\varphi)}\Big)\bigg).
\end{equation*}

The proof of Theorem \ref{thm:inversionf} combined with the convolution theorem
provides an inversion formula for $f$ from $\widehat{R}_\varphi$.
\begin{thm}[Inversion formula for $\widehat{R}_\varphi f$]
Let $f\in C_c^{\infty}(\R^2)$. Then
\begin{equation*}
f(x)=\frac{1}{4\pi}R^* (\widehat{I}^{-1} \widehat{R}_\varphi f)(x).
\end{equation*}
\begin{proof}
By Theorem \ref{modified:form2} and the convolution theorem,
\begin{equation}\label{conv}
F_1(\widehat{R}_\varphi f(\theta,\cdot))(s)
= F_1(Rf(\theta,\cdot))(s) F_1 (\varphi)(s).
\end{equation}
Since $F_1 (\varphi)(s)\neq 0$, by Theorem \ref{thm:pojection-slice},
it follows that
\begin{equation}
F_2f(s\omega) = \frac{F_1(\widehat{R}_\varphi f(\omega,\cdot))(s)}{\sqrt{2\pi}F_1 (\varphi)(s)}. \label{eq:mod-proj-slice}
\end{equation}
Applying the Fourier inversion formula and Proposition \ref{modified:form2} , one shows that
\begin{align*}
f(x) &= \frac{1}{2\pi} \int_{\R^2} F_2f(\xi) e^{i<x,\xi>} d\xi\\
&= \frac{1}{2\pi}\frac{1}{2} \int_{\theta=0}^{2\pi}
\int_{s=-\infty}^{\infty} F_2f(s\omega) e^{i<x,s\omega>} |s| ds d\theta\\
&= \frac{1}{4\pi} \int_{\theta=0}^{2\pi}\int_{s=-\infty}^{\infty}
\widehat{I}^{-1}\widehat{R}_\varphi f(\theta,\langle x, \omega \rangle)d\theta\\
&=\frac{1}{4\pi} R^* (\widehat{I}^{-1} \widehat{R}_\varphi f)(x).
\end{align*}
\end{proof}
\end{thm}

%
\section{Recovering moments from the modified Radon transform}


Let $\omega=(\cos{\theta},\sin{\theta})$ denote a unit direction vector and $x=(x_1, x_2)$ a vector in $\R^2$.
Suppose that $f:\R^2 \rightarrow \R$ is an nonnegative and square-integrable function.
Using the definition of Radon transform, we have 
\begin{equation}\label{range}
\int_{-\infty}^{\infty} Rf(\theta,p)p^k~dp = \int_{\R^2} f(x)\langle \omega, x \rangle^k dx \quad \text{for each }k\in \N_0.
\end{equation}
Assume that the support of $f$ is contained within $I^2$ with $I=[0,1]$.
Then the identity \eqref{range} can be expressed as
\begin{equation}\label{system1}
b^{(k)}(\theta) := \int_{-1}^{1} Rf(\theta,p) p^k~dp =\sum_{j=0}^k
C(k,j)
\big(\cos^{j}{\theta}\sin^{k-j}{\theta}\big)\gamma_{j,k-j},
\end{equation}
where 
\begin{equation*}
C(k,j)=\dfrac{k!}{j!(k-j)!}
\quad \text{and} \quad
\gamma_{\alpha_1,\alpha_2}=\int_{I^2} x_1^{\alpha_1} x_2^{\alpha_2} f(x)dx~
\textup{ for all } \alpha_1,\alpha_2\in \N_0.
\end{equation*}

Let $0< \theta_0 < \theta_{1} < \cdots < \theta_{k} < \pi$ be distinct angles.
Then, system \eqref{system1} can be written in matrix form as follows:
\begin{equation}\label{axb2}
\mathbf{A}^{(k)}\mathbf{x}^{(k)}=\mathbf{b}^{(k)},
\end{equation}
where
\begin{equation*}
\mathbf{A}^{(k)}=
\begin{pmatrix}
C(k,0)\cos^0{\theta_0}\sin^k{\theta_0} & C(k,1)\cos^1{\theta_0}\sin^{k-1}{\theta_0} & \cdots & C(k,k)\cos^k{\theta_0}\sin^0{\theta_0}\\
C(k,0)\cos^0{\theta_1}\sin^k{\theta_1} & C(k,1)\cos^1{\theta_1}\sin^{k-1}{\theta_1} & \cdots & C(k,k)\cos^k{\theta_1}\sin^0{\theta_1}\\
C(k,0)\cos^0{\theta_2}\sin^k{\theta_2} & C(k,1)\cos^1{\theta_2}\sin^{k-1}{\theta_2} & \cdots & C(k,k)\cos^k{\theta_2}\sin^0{\theta_2}\\
\vdots & \vdots & \ddots & \vdots \\
C(k,0)\cos^0{\theta_k}\sin^k{\theta_k} & C(k,1)\cos^1{\theta_k}\sin^{k-1}{\theta_k} & \cdots & C(k,k)\cos^k{\theta_k}\sin^0{\theta_k}
\end{pmatrix},
\end{equation*}

\begin{equation*}
\mathbf{x}^{(k)}=
\begin{pmatrix}
\gamma_{0,k}\\
\gamma_{1,k-1}\\
\gamma_{2,k-2}\\
\vdots\\
\gamma_{k,0}\\
\end{pmatrix},
\mathbf{b}^{(k)}=
\begin{pmatrix}
b^{(k)}(\theta_0)\\
b^{(k)}(\theta_1)\\
b^{(k)}(\theta_2)\\
\vdots\\
b^{(k)}(\theta_k)\\
\end{pmatrix}.
\end{equation*}
The determinant of the matrix $A^{(k)}$ can be expressed as:
\begin{equation*}
\det(A^{(k)})=\det(V^{(k)}) \prod_{j=1}^k \big( C(k,j) \sin^k{\theta_j} \big) ,
\end{equation*}
where $V^{(k)}=[\cot^{j-1}{\theta_i}]_{1\leq i,j \leq k+1}$ is a Vandermonde matrix.
Using the Vandermonde determinant formula, it is easy to show $\det(A^{(k)})$ is positive, implying
the system (\ref{axb2}) has a unique solution $\mathbf{x}^{(k)}$.
Note that the matrix $A^{(k)}$ is positive definite since its leading principal minors are all positive.

Now consider modified Radon transform with noisy data.
\begin{thm}\label{modified:range1}
Assume that the support of $f$ is contained within $[0,1]^2$.
Then for each $k\in \N_0$,
\begin{equation}\label{modified:range2}
\hat{b}^{(k)}(\theta):=\int_{-\infty}^{\infty} \widehat{R}_\varphi f(\theta,p)p^kdp
= \sum_{j=0}^k C(k,j) c_{j}
b^{(k-j)}(\theta)
\end{equation}
where ${\displaystyle c_{j}=\int_{-\infty}^{\infty} \varphi(\tau)(-\tau)^j d\tau}$
for each $j\in \N_0$. 
That is, the value $\hat{b}^{(k)}(\theta)$ is the linear combination of $b^{(0)}(\theta)$, $b^{(1)}(\theta)$, $\ldots$, $b^{(k)}(\theta)$ for any $\theta$.
\begin{proof}
By Lemma \ref{modified:form1}, it follows that
\begin{align*} 
\int_{p=-\infty}^{\infty} \widehat{R}_\varphi  f(\theta,p)p^k~dp
&= \int_{p=-\infty}^{\infty}
\Bigg(\int_{\tau=-\infty}^\infty Rf(\theta,p+\tau) \varphi (\tau) d\tau \Bigg) p^k dp\\
&= \int_{p=-\infty}^{\infty}  \int_{\tau=-\infty}^\infty
\Bigg(\int_{\langle x,\omega\rangle =p+\tau} f(x)dm(x) \Bigg) \varphi (\tau) d\tau~ p^k dp\\
&= \int_{\tau=-\infty}^\infty \varphi (\tau) \Bigg[ \int_{p=-\infty}^{\infty} p^k
\Bigg(\int_{\langle x,\omega\rangle =p+\tau} f(x)dm(x) \Bigg) dp \Bigg] d\tau\\
&= \int_{\tau=-\infty}^\infty \varphi (\tau)
\Bigg(\int_{\R^2} (\langle x,\omega\rangle -\tau)^k f(x)dx \Bigg) d\tau.\\
\end{align*}
Using the following polynomial expansion
\begin{align*}
(\langle x,\omega\rangle-\tau)^k = \sum_{|\alpha|=k} \frac{k!}{\alpha!}
(x_1\cos{\theta})^{\alpha_1}(x_2\sin{\theta})^{\alpha_2}(-\tau)^{\alpha_3},
\end{align*}
where $\alpha=(\alpha_1,\alpha_2,\alpha_3)$ is a multi-index
with $|\alpha|:= \alpha_1 + \alpha_2 + \alpha_3$ and $\alpha!:=\alpha_1!\alpha_2!\alpha_3!$,
it follows that
\begin{align*}
\int_{-\infty}^{\infty} \widehat{R}_\varphi f(\theta,p)p^kdp
&= \sum_{|\alpha|=k}\frac{k!}{\alpha!}c_{\alpha_3}
\gamma_{\alpha_1,\alpha_2} \cos^{\alpha_1}{\theta} \sin^{\alpha_2}{\theta}\\
&= \sum_{\alpha_3=0}^k C(k,\alpha_3) c_{\alpha_3}
\sum_{j=0}^{k-\alpha_3}
C(k-\alpha_3,j)
\gamma_{j,k-\alpha_3-j} \cos^{j}{\theta} \sin^{k-\alpha_3-j}{\theta}.
\end{align*}
\end{proof}
\end{thm}
Let $0< \theta_1 < \theta_{2} < \cdots < \theta_{k} < \pi$ be distinct angles.
Then the system \eqref{modified:range2} can be written in matrix form as follows:
\begin{equation}\label{bbc}
\mathbf{\widehat{B}}^{(k)} = \mathbf{B}^{(k)}\mathbf{C}^{(k)},
\end{equation}
where
\begin{equation*}
\mathbf{C}^{(k)}=
\begin{pmatrix}
C(k,0)c_0 & 0 & \cdots & 0\\
C(k,1)c_1 &  C(k-1,0)c_0 & \cdots & 0\\
C(k,2)c_2 & C(k-1,1)c_1 & \cdots & 0\\
\vdots & \vdots & \ddots & \vdots \\
C(k,k)c_k & C(k-1,k-1)c_{k-1} & \cdots & C(0,0)c_0\\
\end{pmatrix},
\end{equation*}
\begin{equation*}
\mathbf{B}^{(k)}=
\begin{pmatrix}
b^{(k)}(\theta_0) & b^{(k-1)}(\theta_0) & \cdots & b^{(0)}(\theta_0)\\
b^{(k)}(\theta_1) & b^{(k-1)}(\theta_1) & \cdots & 0\\
b^{(k)}(\theta_2) & b^{(k-1)}(\theta_2) & \cdots & 0\\
\vdots & \vdots & \ddots & \vdots \\
b^{(k)}(\theta_k) & 0 & \cdots & 0
\end{pmatrix},
\end{equation*}
\begin{equation*}
\mathbf{\widehat{B}}^{(k)}=
\begin{pmatrix}
\hat{b}^{(k)}(\theta_0) & \hat{b}^{(k-1)}(\theta_0) & \cdots & \hat{b}^{(0)}(\theta_0)\\
\hat{b}^{(k)}(\theta_1) & \hat{b}^{(k-1)}(\theta_1) & \cdots & 0\\
\hat{b}^{(k)}(\theta_2) & \hat{b}^{(k-1)}(\theta_2) & \cdots & 0\\
\vdots & \vdots & \ddots & \vdots \\
\hat{b}^{(k)}(\theta_k) & 0 & \cdots & 0
\end{pmatrix}.
\end{equation*}

Since the matrix $\mathbf{C}^{(k)}$ is a lower triangular matrix, the determinant of the matrix is $(c_0)^{k+1}>0$.
Then there exists a unique solution $\mathbf{B}^{(k)}$ for the given matrix $\mathbf{\widehat{B}}^{(k)}$.
Thus we have the vectors $\mathbf{b}^{(i)}$ for all $i=0,1,\ldots,k$. 
By the equation \eqref{axb2}
one can find $\mathbf{x}^{(0)}$, $\mathbf{x}^{(1)}$, \ldots, $\mathbf{x}^{(k)}$ such that
$\mathbf{x}^{(i)}=(\mathbf{A}^{(i)})^{-1}\mathbf{b}^{(i)}$ for all $i=0, 1,\ldots, k$. That is,
one has the moments $\{\gamma_{\alpha_1,\alpha_2} \}_{\alpha_1+\alpha_2 \leq k}$.

%
%
\section{Recovery of function via moments with noisy data}
Supposed that a cumulative distribution function (cdf) $F$ is absolutely continuous with respect to the Lebesgue measure with the corresponding density function $f$. 
Suppose that $f$ is an $M$-determinate measurable function with a compact support $[0,1]^2$. 
Let $\{\gamma_{\alpha_1,\alpha_2} \}_{\alpha_1\leq m,\alpha_2 \leq n}$ be a given sequence of moments up to $m+n$.
In \cite{ML13} R. Mnatsaknov and S. Li construct the approximation of $f$ such that
\begin{equation*}
\gamma_{\alpha_1,\alpha_2} = \int\int x_1^{\alpha_1} x_2^{\alpha_2} f(x_1,x_2) dx_1dx_2 
\end{equation*}
for all $\alpha_1, \alpha_2\in\N_0$ with $\alpha_1\leq m,~\alpha_2 \leq n$. Their approximation of $f$ is
\begin{equation*}
\frac{\Gamma(m+2)\Gamma(n+2)}{\Gamma([m x_1] +1) \Gamma([n x_2] +1)}
\sum_{\alpha_1=0}^{m-[mx_1]} \sum_{\alpha_2=0}^{n-[nx_2]}
\frac{(-1)^{\alpha_1+\alpha_2} \gamma_{\alpha_1 + [mx_1],\alpha_2 + [nx_2]} }{\alpha_1!\alpha_2!(m-[mx_1]-\alpha_1)!(n-[nx_2]-\alpha_2)!},
\end{equation*}
denoted by $app(f)$.

\begin{thm}\cite{ML13}\label{approx1}
Assume that the pdf $f$ is continuous on $[0,1]^2$. Then $app(f) $ converges uniformly to $ f $ as $m,n \rightarrow \infty$. Furthermore, for any $\delta>0$
\begin{equation*}
||app(f) - f || \leq \Delta(f,\delta) + \frac{4||f||}{\delta^2(\alpha^* +2)}
+\frac{2||f||}{\delta^4 (m+2)(n+2)},
\end{equation*}
where $||\cdot||$ is the sup-norm, $\alpha^*=\min{(m,n)}$, and $\Delta(f,\delta)$ represents the modulus of continuity of $f$.
\end{thm}

\begin{thm}
Assume that the pdf $f$ is continuous on $[0,1]^2$. Then $app( f \ast \varphi) $ converges uniformly to $ f $ on $ [0,1]^2 $.
\begin{proof}
By Theorem \ref{approx1}
it follows that for any $\delta>0$
\begin{align*}
||f-app(f\ast \varphi)||
& \leq ||f-app(f)||+ ||app(f)-app(f\ast\varphi)||\\
& \leq \Delta(f,\delta) + \frac{4(||f||}{\delta^2(\alpha^* +2)}
+\frac{2(||f||}{\delta^4 (m+2)(n+2)} +
\end{align*}
\end{proof}
\end{thm}





\end{document}